\documentclass[a4paper,11pt]{article}
\usepackage[latin1]{inputenc}
\usepackage{amsfonts}
\usepackage{amsmath}
\usepackage{amssymb}
\usepackage{amsthm}
\usepackage[T1]{fontenc}
\usepackage[dvips]{graphicx}
\usepackage{color}
\usepackage{epsfig}
\theoremstyle{plain}
\newtheorem{Thm}{Theorem}
\newtheorem{Lem}{Lemma}
\newtheorem{Def}{Definition}

\newtheorem*{Rem*}{\textsc{Remark}}
\newtheorem*{Def*}{Definition}
\newtheorem*{Prob*}{Direct Monodromy Problem}
\newtheorem*{Lem*}{\textsc{Lemma}}
\newtheorem*{Cor*}{\textsc{Corollary}}
\newtheorem*{Con*}{\textsc{Conjecture}}
\newcommand{\e}{\varepsilon}
\newcommand{\bb}[1]{\mathbb{ #1 }}
\newcommand{\av}[1]{\left| #1 \right|}

\title{The Direct Monodromy Problem of Painleve-I}
\author{Davide Masoero \thanks{E-mail address: masoero@sissa.it}\\ SISSA, Trieste
and Grupo de F\'isica Matem\'atica, Lisboa}
\date{}
\begin{document}

\maketitle

\abstract{The Painleve first equation can be represented as the equation of isomonodromic deformation
of a Schr\"odinger equation with a cubic potential. We introduce a new algorithm for computing the direct
monodromy problem for this Schr\"odinger equation. The algorithm is based on the geometric theory of 
Schr\"odinger equation due to Nevanlinna}

\section{Introduction}

Painlev\'e equations form the core of what may be called "modern
special function theory". Indeed, since the 1970s' pioneering works of Ablowitz and Segur, and McCoy, Tracy and Wu
the Painlev\'e functions have been playing the same role in nonlinear mathematical physics that the classical special functions, such as Airy functions, Bessel functions, etc., are playing in linear physics.

For example, it has been recently conjectured by Dubrovin \cite{dubrovin08} (and proven in a particular case by Claeys and Grava \cite{tamtom}) that Painlev\'e equations play a big role also in the theory of nonlinear waves and dispersive equations.

In this context, the first Painlev\'e equation (P-I)
\begin{equation}\label{eq:PI}
y''= 6y^2 -z \, , \; z \in \mathbb{C} \quad ,
\end{equation}
is of special importance.

Indeed recently
\cite{dubrovin} Dubrovin, Grava and Klein
discovered that a special solution of P-I -called int\'egrale tritronqu\'ee- provides
the universal correction to the semiclassical limit of solutions to the focusing nonlinear
Schr\"odinger equation.

The key mathematical fact about the Painlev\'e equations is their Lax-integrability.
This allows to apply to the study of the Painlev\'e functions a powerful Isomonodromy - Riemann-Hilbert method.
Using this method, a great number of analytical, and especially asymptotic results have been obtained during the last
two decades- for the general theory see the monumental book \cite{fokas}, for what concerns P-I see
our papers \cite{piwkb} \cite{piwkb2}.

It is clear that we need efficient and reliable algorithms for computing solutions
of the Painlev\'e equations and that these algorithms should take into account the
integrability of Painlev\'e equations. To this aim S. Olver has recently build an algorithm
to solve the Riemann-Hilbert problem -or inverse monodromy problem-
associated to Painlev\'e-II \cite{sheehan10}.

In the present paper we introduce a new simple algorithm
for solving the \textit{direct monodromy problem} associated to P-I;
given the Cauchy data $y,y',z$ of a given solution of P-I we can compute
its monodromy data -or equivalently, the Riemann-Hilbert problem associated to it.

Given our algorithm, it would be possible to compute the inverse monodromy problem
for Painlev\'e-I simply using an appropriate Newton's method.
Due to the simplicity of our algorithm, this procedure would be probably more efficient than a numerical solution
of the corresponding Riemann-Hilbert problem (which, to the best of our knowledge, has not yet been given).
However, this is a matter of future research.

Below we briefly introduce the direct monodromy problem associated to P-I.

\subsection{The Perturbed Cubic Oscillator}
Consider the following Schr\"odinger equation with a cubic potential (plus a fuchsian singularity)
\begin{eqnarray}\label{eq:perturbedschr} 
 \frac{d^2\psi(\lambda)}{d\lambda^2} \!\! &=& \!\! Q(\lambda;y,y',z) \psi(\lambda)\; , \\ \nonumber
 Q(\lambda;y,y',z) \!\! &=& \!\! 4 \lambda^3 - 2 \lambda z + 2 z y - 4 y^3 + y'^2+  \frac{y'}{\lambda - y} 
+\frac{3}{4(\lambda -y)^2} \, .
\end{eqnarray}
We call such equation the \textit{perturbed (cubic) oscillator}.

We define the Stokes Sector $S_k$ as
\begin{equation}
 S_k=\left\lbrace \lambda : \av{ \arg\lambda - \frac{2 \pi k}{5} } < \frac{\pi}{5} \right\rbrace
 \, , k \in \bb{Z}_{5} \; .
\end{equation}
Here, and for the rest of the paper, $\bb{Z}_5$ is the group of the integers modulo five. We will often choose
as representatives of $\bb{Z}_{5}$ the numbers $-2,-1,0,1,2$.

For any Stokes sector, there is a unique (up to a multiplicative constant) solution of the perturbed oscillator
that decays exponentially inside $S_k$. We call such solution the \textit{k-th subdominant solution} and let 
$\psi_k(\lambda;y,y',z)$ denote it \footnote{Equation (\ref{eq:perturbedschr}) has a fuchsian singularity
at $\lambda=y$. Hence, as we explain in Section 2, any solution is two-valued having a square root singularity at $\lambda=y$.
We do not discuss this fact here: the reader can find the rigorous definition of subdominant solution in the Appendix.}.

The asymptotic behaviour of $\psi_k$ is known explicitely in a bigger sector
of the complex plane, namely $S_{k-1}\cup\overline{S_k}\cup S_{k+1}$:
\begin{equation}\label{eq:intwkb}
\lim_{\substack{ \lambda \to \infty \\ \av{\arg{\lambda} - \frac{2 \pi k}{5}} < \frac{3 \pi}{5} -\e}}
\frac{\psi_k(\lambda;y,y',z)}{\lambda^{-\frac{3}{4}} \exp\left\lbrace-\frac{4}{5} \lambda^{\frac{5}{2}} +
\frac{z}{2}\lambda^{\frac{1}{2}}\right\rbrace} \to 1 , \; \forall \e >0 \, .
\end{equation}
Here the branch of $\lambda^{\frac{1}{2}}$ is chosen such that $\psi_k$ is exponentially small in $S_k$.

Since $\psi_{k-1}$ grows exponentially in $S_k$, then $\psi_{k-1}$ and $\psi_{k}$ are linearly independent.
Then $\left\lbrace \psi_{k-1},\psi_{k} \right\rbrace$ is a basis of solutions, whose asymptotic behaviours is known in $S_{k-1}\cup S_{k}$.

Fixed $k^* \in \bb{Z}_5$, we know the asymptotic behaviour of $\left\lbrace \psi_{k^*-1},\psi_{k^*} \right\rbrace$
only in $S_{k^*-1}\cup S_{k^*}$.
If we want to know the asymptotic behaviours of this basis in all the complex plane, it is sufficient to
know the linear transformation from basis $\left\lbrace \psi_{k-1},\psi_{k} \right\rbrace$ to basis
$\left\lbrace \psi_{k},\psi_{k+1} \right\rbrace$ for any $k \in \bb{Z}_5$.

From the asymptotic behaviours, it follows that these changes of basis
are triangular matrices: for any $k$, $\psi_{k-1}=\psi_{k+1}+\sigma_k\psi_k$
for some complex number $\sigma_k$, called \textit{Stokes multiplier}.
The quintuplet of Stokes multipliers $\sigma_k, k \in \bb{Z}_5$ is called the monodromy data of the perturbed oscillator.

We can now define the direct monodromy problem.
\begin{Prob*}
Fixed $y,y',z$ compute the Stokes multipliers $\sigma_k(y,y',z)$
of the perturbed oscillator equation (\ref{eq:perturbedschr}).
\end{Prob*}
Our Algorithm gives a numerical solution
of this problem.

\subsection{P-I: Isomonodromy Approach}

What is the relation of the direct monodromy problem with P-I?
\textbf{Painlev\'e-I is the equation of isomonodromic deformation of the perturbed oscillator}.
Indeed the following Theorem holds.
\begin{Thm}\label{thm:isomonodromy}
Let the parameters $y=y(z), y'=\frac{dy(z)}{dz}$ of the potential $Q(\lambda;y,y',z)$
be functions of $z$; then $y(z)$
solves P-I if and only if the Stokes multipliers of the perturbed oscillator do not depend on $z$.
\begin{proof}
See for example \cite{piwkb}.
\end{proof}
\end{Thm}

Fix a solution $y(z)$ of P-I. If $z$ is a pole of $y(z)$ then equation (\ref{eq:perturbedschr}) is not well-defined. However,
recently the author \cite{piwkb} (see also \cite{chudnovsky94}) showed that
this difficulty can be overcome.

Let $a \in \bb{C}$ be a pole of $y(z)$ then
$y(z)$ has the following Laurent expansion around $a$
$$y(z)=\frac{1}{(z-a)^2}+\frac{a(z-a)^2}{10}+ \frac{(z-a)^3}{6}+b(z-a)^4+O((z-a)^5) \; .$$

Inserting the Laurent expansion into equation (\ref{eq:perturbedschr}) and taking the the limit $z \to a$,
equation (\ref{eq:perturbedschr})
becomes the following Schr\"odinger equation with a cubic potential
\begin{equation}\label{eq:schr}
\frac{d^2\psi(\lambda)}{d\lambda^2}= V(\lambda;2a,28b) \psi(\lambda)\; , \quad
V(\lambda;a,b)=4 \lambda^3 - a \lambda -b \;.
\end{equation}
Moreover the Stokes multipliers of the equation (\ref{eq:schr}) coincide with
the Stokes multipliers of equation (\ref{eq:perturbedschr}).

We call equation (\ref{eq:schr}) the cubic oscillator. For the rest of the paper
we consider it as a particular case of the perturbed oscillator (\ref{eq:perturbedschr}).

\subsection{Stokes Multipliers and Asymptotic Values}
The algorithm is based on the formula (\ref{eq:intsigmaR})
below, that we discovered in \cite{dtba}.

Consider the following Schwarzian equation
\begin{equation}\label{eq:intschw}
\left\lbrace f(\lambda) , \lambda \right\rbrace = -2 Q(\lambda;y,y',z) \, .
\end{equation}
Here $\left\lbrace f(\lambda) , \lambda \right\rbrace=\frac{f'''(\lambda)}{f'(\lambda)} -
\frac{3}{2}\left(\frac{f''(\lambda)}{f'(\lambda)}\right)^2 $ is the Schwarzian derivative.

For every solution of the Schwarzian equation (\ref{eq:intschw}) the following limit exists
\begin{equation*}
w_k(f)=\lim_{\lambda\ \to \infty \,, \lambda \in S_k }f(\lambda) \in \mathbb{C}
\cup \infty \, ,
\end{equation*}
provided the limit is taken along a curve non-tangential to the boundary of $S_k$.

In Section 2, we will prove that the following formula holds for any solution of the Schwarzian
equation (\ref{eq:intschw})
\begin{equation}\label{eq:intsigmaR}
\sigma_k(y,y',z)=i\left(w_{1+k}(f),w_{-2+k}(f); w_{-1+k}(f),w_{2+k}(f) \right) \;.
\end{equation}
Here $(a,b;c,d)=\frac{(a-c)(b-d)}{(a-d)(b-c)}$ is the cross ratio of four points on the sphere.

The paper is organized as follows.
In Section 2 we introduce the Nevanlinna's theory of the cubic oscillator and the Schwarzian differential equation (\ref{eq:schw}).
Then we prove formula (\ref{eq:intsigmaR}) formula for computing the Stokes multipliers from any solution of the Schwarzian differential equation.
Section 3 is devoted to the description of the Algorithm.
In Section 4 we test our algorithm against the WKB prediction and the Deformed TBA equations.
For convenience of the reader, we explain the basic theory of cubic oscillators (Stokes sectors, Stokes multipliers,
subdominant solutions, etc ...) in the Appendix.

\paragraph{Acknowledgments}
I am indebted to my advisor Prof. B. Dubrovin who constantly gave me
suggestions and advice. This work began in May 2010 during the workshop "Numerical solution of the Painlev\'e equations" at ICMS, Edinburgh  and was finished
in June 2010 while I was a guest of Prof. Y. Takei at RIMS, Kyoto.
I thank Prof. Takei and all the participants to the workshop for the stimulating discussions. 
This work is partially supported by the Italian Ministry of University and Research
(MIUR) grant PRIN 2008 "Geometric methods in the theory of nonlinear waves and their applications".

\section{Schwarzian Differential Equation}
As we mentioned in the Introduction, our Algorithm is based on formula (\ref{eq:intsigmaR}) that allows one
to compute Stokes multipliers from any solution of the Schwarzian differential equation (\ref{eq:intschw}).

This formula, proven in Theorem \ref{thm:dtba} below,
has its roots in the geometric theory of the Schr\"odinger equation,
which was developed by Nevanlinna in the 1930s' \cite{nevanlinna32}.
The author learned such a beautiful theory from the remarkable paper of Eremenko and Gabrielov \cite{eremenko}.
In this section we follow quite closely \cite{eremenko} as well as author's recent paper \cite{dtba}.

\begin{Rem*}
Equation (\ref{eq:perturbedschr}) has a fuchsian singularity at the pole $\lambda=y$ of the potential
$Q(\lambda;y,y',z)$. However this is an \textit{apparent singularity} \cite{piwkb}: the monodromy around
the singularity of any solution of (\ref{eq:perturbedschr}) is $-1$.
As a consequence, the ratio of two solutions of (\ref{eq:perturbedschr}) is a meromorphic function.
\end{Rem*}

The main geometric object of Nevanlinna's theory is the Schwarzian derivative of a
(non constant) meromorphic function $f(\lambda)$

\begin{equation}\label{def:schwarzian}
\left\lbrace f(\lambda),\lambda \right\rbrace =\frac{f'''(\lambda)}{f'(\lambda)} -
\frac{3}{2}\left(\frac{f''(\lambda)}{f'(\lambda)}\right)^2 \; .
\end{equation}

The Schwarzian derivative is strictly related to the Schr\"odinger equation (\ref{eq:schr}). Indeed, the following Lemma is
true.

\begin{Lem}\label{lem:schsch}
The (non constant) meromorphic function $f:\bb{C} \to \overline{\bb{C}}$ solves the Schwarzian differential equation
\begin{equation}\label{eq:schw}
\left\lbrace f(\lambda) , \lambda \right\rbrace = -2 Q(\lambda;y,y',z) \, .
\end{equation}
iff $f(\lambda)=\frac{\phi(\lambda)}{\chi(\lambda)}$ where $\phi(\lambda)$ and $\chi(\lambda)$
are two linearly independent solutions of the Schr\"odinger equation (\ref{eq:perturbedschr}). Hence,
the first derivative of any (non constant) solution of (\ref{eq:schw}) vanishes only at the pole $\lambda=y$ of the potential.
\end{Lem}

We define the Asymptotic Stokes Sector $S_k$ as
\begin{equation}
 S_k=\left\lbrace \lambda : \av{ \arg\lambda - \frac{2 \pi k}{5} } < \frac{\pi}{5} \right\rbrace
 \, , k \in \bb{Z}_{5} \; .
\end{equation}

Every solution of the Schwarzian equation (\ref{eq:schw}) has limit for $\lambda \to \infty$, $\lambda \in S_k$.
More precisely we have the following 

\begin{Lem}[Nevanlinna]\label{lem:wk}
\begin{itemize}
\item[(i)]Let $f(\lambda)=\frac{\phi(\lambda)}{\chi(\lambda)}$ be a solution of (\ref{eq:schw}) then for
all $k \in \bb{Z}_5$ the following limit exists
\begin{equation}\label{eq:wk}
w_k(f)=\lim_{\lambda\ \to \infty \, \lambda \in S_k }f(\lambda) \in \mathbb{C}
\cup \infty \, ,
\end{equation}
provided the limit is taken along a curve non-tangential to the boundary of $S_k$.
\item[(ii)]$w_{k+1}(f) \neq w_{k}(f) \, , \; \forall k \in \bb{Z}_5$.
 \item[(iii)] Let $g(\lambda)=\frac{a f(\lambda) +b}{c f(\lambda)+d}=
 \frac{a\phi(\lambda) + b\chi(\lambda) }{c \phi(\lambda)+d\chi(\lambda)}$, 
$\left(\begin{matrix}
a & b \\ 
c  & d
\end{matrix} \right) \in Gl(2,\bb{C})$. Then
\begin{equation}\label{eq:moebius}
w_k(g)= \frac{a \, w_k(f) + b}{c \, w_k(f)+d} \; .
\end{equation}
\item[(iv)]If the function $f$ is evaluated along a ray contained in $S_k$, the convergence to $w_k(f)$ is super-exponential.
\end{itemize}
\begin{proof}

\begin{itemize}
 \item[(i-iii)]Let $\psi_k$ be the solution of equation (\ref{eq:perturbedschr}) subdominant in $S_k$ and $\psi_{k+1}$ be the one subdominant in $S_{k+1}$. Since they form a basis of solutions, then $f(\lambda)=\frac{\alpha\psi_k(\lambda) + \beta\psi_{k+1}(\lambda) }{\gamma \psi_k(\lambda)+\delta\psi_{k+1}(\lambda)}$, for some $\left(\begin{matrix}
\alpha & \beta \\ 
\gamma  & \delta
\end{matrix} \right) \in Gl(2,\bb{C}) $. Hence $w_k(f)=\frac{\beta}{\delta}$ if $\delta \neq 0$, $w_k(f)=\infty$ if $\delta=0$.
Similarly $w_{k+1}(f)=\frac{\alpha}{\gamma}$. Since $\left(\begin{matrix}
\alpha & \beta \\ 
\gamma  & \delta
\end{matrix} \right) \in Gl(2,\bb{C}) $ then $w_k(f) \neq w_{k+1}(f)$
\item[(iv)]From estimates (\ref{eq:intwkb}) we know that inside $S_k$, 
$$\av{\frac{\psi_k(\lambda;y,y',z)}{\psi_{k+1}(\lambda;y,y',z)}} \sim e^{-Re\left(\frac{8}{5}\lambda^{\frac{5}{2}} -
z\lambda^{\frac{1}{2}}\right)} \;,$$ where the branch of $\lambda^{\frac{1}{2}}$ is chosen such that
the exponential is decaying.
\end{itemize}

\end{proof}
\end{Lem} 

\begin{Def}\label{def:wk}
Let $f(\lambda)$ be a solution of the Schwarzian equation (\ref{eq:schw})
and $w_k(f)$ be defined as in (\ref{eq:wk}). We call $w_k(f)$ the k-th asymptotic value of $f$.
\end{Def}

We noticed in a previous paper \cite{dtba} that the Stokes multipliers of the Schr\"odinger equation
are rational functions of the asymptotic values $w_k(f)$. This relation is the basis of our Algorithm.

\begin{Thm}\label{thm:dtba} \cite{dtba}
Denote $\sigma_k$ the k-th Stokes multiplier of the Schr\"odinger equation (\ref{eq:perturbedschr})
(for its precise definition, see equation (\ref{eq:multipliers}) in the Appendix).
Let $f$ be any solution of the Schwarzian equation (\ref{eq:schw}). Then
\begin{equation}\label{eq:identity}
\sigma_k = i \left(w_{1+k}(f),w_{-2+k}(f); w_{-1+k}(f),w_{2+k}(f) \right) \;, \forall k \in \bb{Z}_5 \;,
\end{equation}
where $(a,b;c,d)=\frac{(a-c)(b-d)}{(a-d)(b-c)}$ is the cross ratio of four point on the sphere.

\begin{proof}
Due to equation (\ref{eq:moebius}) all the asymptotic values
of two different solutions of (\ref{eq:schw}) are related by the same fractional linear transformation.
As it is well-known, the cross ratios of four points of the sphere is invariant if all the points are transformed by the same fractional linear transformation.
Hence the right-hand side of (\ref{eq:identity}) does not depend on the choice of the solution of the Schwarzian equation.

Let $\psi_{k+1}$ be the solution of (\ref{eq:perturbedschr}) subdominant in $S_{k+1}$ and $\psi_{k+2}$ be
the one subdominant in $S_{k+2}$ (see the Appendix for the precise definition). By choosing
$f(\lambda)=\frac{\psi_{k+1}(\lambda)}{\psi_{k+2}(\lambda)}$, one verifies easily that the
identity (\ref{eq:identity}) is satisfied.
\end{proof}

\end{Thm}

\begin{Rem*}
The same construction presented here holds for anharmonic oscillators with polynomial potentials
of any degree. For any degree, there are formulas similar to (\ref{eq:identity}) for expressing
Stokes multipliers in terms of cross ratios of asymptotic values.
The general formula will be given in a subsequent publication.
\end{Rem*}

\subsection{Singularities}

Since the Schwarzian differential equation is linearized  (see Lemma \ref{lem:schsch}) by the Schr\"odinger equation, any solution is a meromorphic function and has an infinite number of poles (for a proof of this fact see \cite{nevanlinna32}
and \cite{elfving}).
The poles, however, are  localized near the boundaries of the Stokes sectors $S_k, k \in \bb{Z}_5$.
Indeed, using the estimates (\ref{eq:intwkb}) one can prove the following

\begin{Lem}\label{lem:poles}
Let $f(\lambda)$ be any solution of the Schwarzian equation (\ref{eq:schw}). Fix $\e >0$ and define
$\tilde{S_k}= \left\lbrace \lambda : \av{ \arg\lambda - \frac{2 \pi k}{5} } \leq \frac{\pi}{5} -\e \right\rbrace
 \, , k \in \bb{Z}_{5} \; .$ Then $f(\lambda)$ has a finite number of poles inside $\tilde{S_k}$. Hence,
there are a finite number of rays inside $\tilde{S_k}$ on which $f(\lambda)$ has a singularity.
\end{Lem}

\section{The Algorithm}
In the previous section we have proved the following remarkable facts
\begin{itemize}
 \item Along any ray contained in the Stokes Sector $S_k$, any solution $f$ to the Schwarzian differential equation (\ref{eq:perturbedschr}) $f$
 converges super-exponentially to the asymptotic value $w_k(f)$. See Lemma \ref{lem:wk} (iv).
 \item The Stokes multipliers of the Schr\"odinger equation (\ref{eq:schw}) are cross-ratios of the
 asymptotic values $w_k(f)$. See equation (\ref{eq:identity}).
 \item Inside any closed subsector of $S_k$, $f$ has a finite number of poles. See Lemma \ref{lem:poles}.
\end{itemize}

Hence the Simple Algorithm for Computing Stokes Multipliers goes as follows:

\begin{enumerate}
\item Set k=-2.
\item Fix arbitrary Cauchy data of $f$: $f(\lambda^*),f'(\lambda^*),f''(\lambda^*)$, with the conditions $\lambda^* \neq y$, $f'(\lambda^*) \neq 0$.
\item Choose an angle $\alpha$ inside $S_k$, such that the singular point $\lambda=y$ does not belong to the
corresponding ray, i.e. $\alpha \neq \arg y$. Define $t: \bb{R}^+\cup 0 \to \bb{C}$, $t(x)=f(e^{i\alpha}x +\lambda^*)$.
The function $t$ satisfies the following Cauchy problem
\begin{eqnarray}\label{eq:cauchy}
\left\lbrace \begin{aligned}
   \left\lbrace t(x),x \right\rbrace = e^{2 i\alpha} Q(e^{i\alpha}x + \lambda^*;y,y',z), \hspace{1.5cm}           \\
 t(0)=f(\lambda^*), t'(0)=e^{i\alpha}f'(\lambda^*), \, t''(0)=e^{2 i\alpha}f''(\lambda^*) \; . 
             \end{aligned}
\right.  
\end{eqnarray} 
\item Integrate equation (\ref{eq:cauchy}) either directly \footnote{Integrating equation (\ref{eq:cauchy}) directly, one can hit a singularity $x*$ of $y$. To continue the solution past the pole, starting from $x^*-\e$ one can integrate the function $\tilde{y}=\frac{1}{y}$, which satisfies the same Schwarzian differential equation.} or by linearization (see Remark below), and compute $w_k(f)$ with the desired accuracy and precision.
\item If $k <2, k++$, return to point 3.
\item Compute $\sigma_l$ using formula (\ref{eq:identity}) for all $l \in \bb{Z}_5$.
\end{enumerate}

\begin{Rem*}
As was shown in Lemma \ref{lem:schsch}, any solution $f$ of the Schwarzian equation is the ratio of two solutions of
the Schr\"odinger equation. Hence, one can solve the nonlinear Cauchy problem (\ref{eq:cauchy}) by solving
two linear Cauchy problems.

Whether the linearization is more efficient than the direct integration of (\ref{eq:cauchy}) will not be investigated in the present paper.
\end{Rem*}

\section{A Test}
We have implemented our algorithm using MATHEMATICA's ODE solver NDSOLVE integrating equation (\ref{eq:cauchy})
with steps of length 0.1. We decided the integrator to stop at step $n$ if
\begin{eqnarray*}
\av{t(0.1 n) - t(0.1(n-1))} &<& 10^{-13} \, \mbox{ and } \\
\av{\frac{t(0.1 n) - t(0.1(n-1))}{t(0.1 n)}} &<& 10^{-13} \; .
\end{eqnarray*}

To test our algorithm we computed the Stokes multiplier $\sigma_0(b)$ of the equation

\begin{equation}\label{eq:simpleschr}
\frac{d^2\psi(\lambda)}{d\lambda^2}= (4 \lambda^3 -b) \psi(\lambda) \; .
\end{equation}

According to the WKB analysis (see \cite{sibuya75}, \cite{piwkb}) the Stokes multiplier $\sigma_0(b)$ has the following asymptotics

\begin{equation}\label{eq:wkbprediction}
\sigma_0(b) \sim  \left\lbrace 
             \begin{aligned}
               -i e^{\frac{\sqrt{\frac{\pi}{3}} \Gamma(1/3)}{2^{2/3} \Gamma(11/6)} b^{\frac{5}{6}}} \hspace{2cm} , \; \mbox{ if } \, b>0  \\
 -2 i e^{-\frac{\sqrt{\pi}  \Gamma(1/3)}{ 2^{\frac{5}{3}} \Gamma(11/6)} (-b)^{\frac{5}{6}}} 
\cos\left(\frac{\sqrt{\frac{\pi}{3}} \Gamma(1/3)}{2^{5/3} \Gamma(11/6)} (-b)^{\frac{5}{6}} \right) ,  \mbox{ if } b <0 \; . 
             \end{aligned}
\right. 
\end{equation}

Our computations (see Figure 1 and 2 below) shows clearly that the WKB approximation is very efficient also for small value of the parameter $b$.

\begin{figure}[htbp]
\begin{center}
\includegraphics[width=10cm]{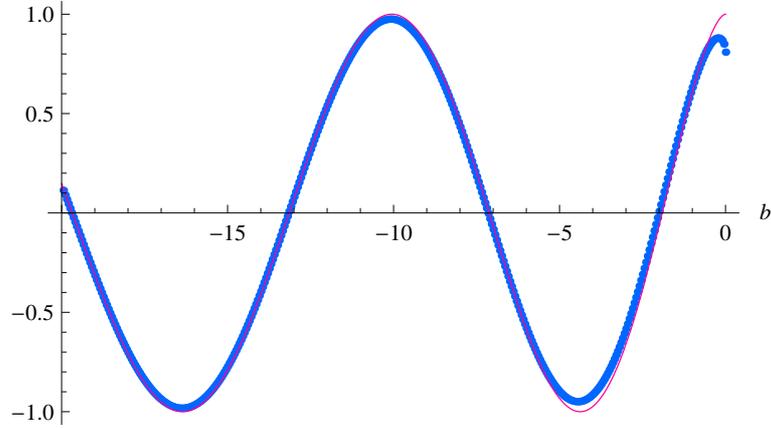}
\end{center}
\label{fig:negativeb}
\caption{ \footnotesize{Thick dotted line: the rescaled Stokes multiplier $\frac{i}{2}e^{\frac{\sqrt{\pi}  \Gamma(1/3)}{ 2^{\frac{5}{3}} \Gamma(11/6)} (-b)^{\frac{5}{6}}} \sigma_0(b)$ evaluated with our algorithm; thin continuous line: $\cos\left(\frac{\sqrt{\frac{\pi}{3}} \Gamma(1/3)}{2^{5/3} \Gamma(11/6)} (-b)^{\frac{5}{6}} \right)$, i.e. the WKB prediction for the rescaled Stokes multiplier.}}
\end{figure}

\begin{figure}[htbp]
\begin{center}
\includegraphics[width=10cm]{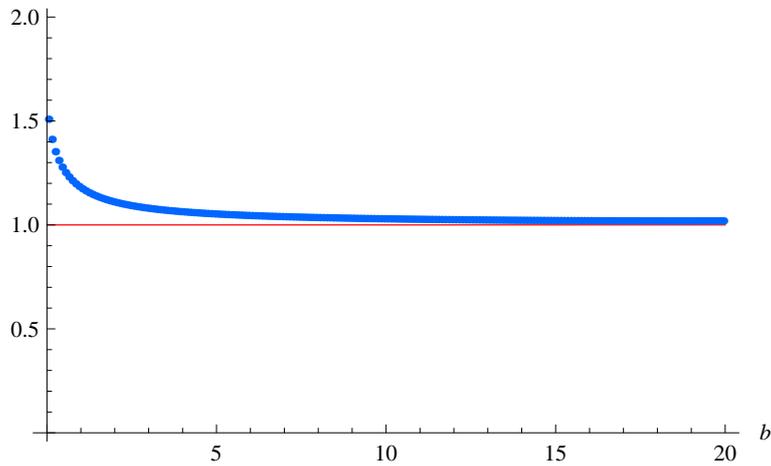}
\end{center}
\label{fig:positiveb}
\caption{ \footnotesize{Thick dotted line: the rescaled Stokes multiplier $e^{-\frac{\sqrt{\frac{\pi}{3}} \Gamma(1/3)}{2^{2/3} \Gamma(11/6)} b^{\frac{5}{6}}} \sigma_0(b)$ evaluated with our algorithm; thin continuous line: 1, the WKB prediction for the rescaled Stokes multiplier.}}
\end{figure}

We also tested our results against the numerical solution (due to A. Moro and the author)
of the Deformed Thermodynamic Bethe Ansatz equations (Deformed TBA), which has been
recently introduced by the author \cite{dtba}, developing the seminal work of Dorey and Tateo \cite{doreytateo}.
The Deformed TBA equations are a set of nonlinear integral equations which describe the exact correction to the WKB asymptotics.
The numerical solution of the Deformed TBA equations enable to a-priori set the absolute error in the evaluation of the Stokes multiplier $\sigma_0(b)$
rescaled with respect to the exponential factor, shown in (\ref{eq:wkbprediction}).
Hence,  in the range of $ -20 \leq b \leq 20$ we could verify that we had computed the rescaled $\sigma_0(b)$ with an absolute error less than $10^{-8}$.

\section{Appendix}

The reader expert in anharmonic oscillators theory will skip this Appendix; for her, it will be enough to know
that we denote  $\sigma_k(y,y',z)$  the $k-th$ Stokes multipliers of equation (\ref{eq:perturbedschr}).
Here we review briefly the standard way, i.e. by means of Stokes multipliers,
of introducing the monodromy problem
for equation (\ref{eq:perturbedschr}). All the statements of this section are proved in
Appendix A of author's paper \cite{piwkb} and in Sibuya's book \cite{sibuya75} .

\begin{Lem}\label{lem:wkb}
Fix $k \in \mathbb{Z}_5 = \left\lbrace -2, \dots , 2 \right\rbrace$ and define a cut in the $\bb{C}$ plane
connecting $\lambda=y$ with infinity such that its points eventually do not belong to
$S_{k-1} \cup \overline{S_k} \cup S_{k+1}$. Choose the branch of $\lambda^{\frac{1}{2}}$ by requiring
$$\lim_{\substack{\lambda \to \infty \\ \arg{\lambda}=
\frac{2 \pi k}{5}}} {\rm Re}\lambda^{\frac{5}{2}} = + \infty \, ,$$
while choose arbitrarily one of the branch of $\lambda^{\frac{1}{4}}$.
Then there exists a unique solution $\psi_k(\lambda;a,b)$ of equation (\ref{eq:schr})
such that
\begin{equation}\label{eq:asym}
\lim_{\substack{ \lambda \to \infty \\ \av{\arg{\lambda} - \frac{2 \pi k}{5}} < \frac{3 \pi}{5} -\e}}
\frac{\psi_k(\lambda;y,y',z)}{\lambda^{-\frac{3}{4}} e^{-\frac{4}{5} \lambda^{\frac{5}{2}} +
\frac{z}{2}\lambda^{\frac{1}{2}}}} \to 1 , \; \forall \e >0 \, .
\end{equation}

\end{Lem}

\begin{Def*}
We denote $\psi_k$ the k-th subdominant solution or the solution subdominant in the k-th sector.
\end{Def*}

From the asymptotics (\ref{eq:asym}), it follows that $\psi_k$ and $\psi_{k+1}$
are linearly independent. If one fixes the same branch of $\lambda^{\frac{1}{4}}$
in the asymptotics (\ref{eq:asym}) of $\psi_{k-1},\psi_{k},\psi_{k+1}$ then the following equations hold true
\begin{eqnarray}\nonumber
  \frac{\psi_{k-1}(\lambda;y,y',z)}{\psi_k(\lambda;y,y',z)} &=& \frac{\psi_{k+1}(\lambda;y,y',z)}{\psi_k(\lambda;y,y',z)}+\sigma_k(y,y',z) \;, \\ \label{eq:multipliers} \\ \nonumber
 - i \sigma_{k+3}(y,y',z)   &=& 1+\sigma_k(y,y',z)\sigma_{k+1}(y,y',z)\, , \; \forall k \in \bb{Z}_5 \; .
\end{eqnarray}

\begin{Def*}
The entire functions $\sigma_k(y,y',z)$ are called Stokes multipliers.
The quintuplet of Stokes multipliers $\sigma_k(y,y',z), k \in \bb{Z}_5$ is called the monodromy data
of equation (\ref{eq:schr}).
\end{Def*}


\end{document}